\documentclass[reqno]{amsart}

\usepackage{dsfont,colonequals,mathrsfs}
\usepackage{amsfonts,amssymb}
\usepackage{enumerate}
\usepackage{units,url}

\title[Approximation of the Robin Semigroup]
	{Approximation of the Semigroup generated by the Robin Laplacian
		in terms of the Gaussian Semigroup}
\author[R.\ Nittka]{Robin Nittka}
\address{University of Ulm, Institute of Applied Analysis, Germany}
\thanks{The author thanks the Graduate School \textsc{Mathematical Analysis of Evolution, Information and Complexity}
		for their support during the work on this article}
\email{robin.nittka@uni-ulm.de}
\date{\today}


\theoremstyle{plain}
\newtheorem{theorem}{Theorem}
\newtheorem{corollary}[theorem]{Corollary}
\newtheorem{proposition}[theorem]{Proposition}
\newtheorem{lemma}[theorem]{Lemma}

\theoremstyle{definition}
\newtheorem{remark}[theorem]{Remark}

\newtheorem{definition}[theorem]{Definition}

\newcommand{\coloneqq}{\colonequals}
\renewcommand{\mid}{\middle|}
\newcommand{\dx}{\mathrm{d}}
\newcommand{\setone}{\mathds{1}}
\newcommand{\eps}{\varepsilon}
\renewcommand{\rho}{\varrho}
\renewcommand{\phi}{\varphi}
\newcommand{\laplace}{\Delta}
\newcommand{\nlaplace}{{\laplace\!}_N}
\newcommand{\rlaplace}{{\laplace\!}_R}
\newcommand{\dlaplace}{{\laplace\!}_D}
\newcommand{\Claplace}{{\laplace}_0}
\newcommand{\Llaplace}{{\laplace}_2}

\newcommand{\Id}{I}

\DeclareMathOperator{\e}{e}

\begin{document}

\begin{abstract}
	For smooth bounded open sets in euclidean space, we construct corresponding
	contractive linear extension operators for the space of continuous functions
	which preserve regularity of functions in the domain of the Robin Laplacian.
	We also prove a Trotter-like approximation
	for the semigroup generated by the Laplacian subject to
	Robin boundary conditions in terms of these extension operators.
	The limiting case of Dirichlet boundary conditions is treated
	separately.
\end{abstract}

\maketitle

\section{Introduction}
Let $\Omega \subset \mathds{R}^N$ be a smooth bounded open set. Here and in the following,
``smooth'' means ``of class $\mathrm{C}^\infty$'', although the main results remain true under slightly
milder regularity assumptions.
On such a set we consider the (autonomous, homogeneous) diffusion equation
\begin{equation}\label{diffeq}
	\left\{
		\begin{aligned}
			u_t & = \laplace u, && \text{on } (0,\infty) \times \Omega, \\
			\frac{\partial u}{\partial\nu}(t,z) & = -\beta(z) u(t,z), && \text{for } t > 0 \text{ and } z \in \partial\Omega, \\
			u(0,x) & = u_0(x), && \text{for } x \in \Omega,
		\end{aligned}
	\right.
\end{equation}
subject to Robin boundary conditions.
Here $u_0 \in \mathrm{C}(\overline{\Omega})$ is an arbitrary initial function,
$\beta$ is a non-negative smooth function on $\partial\Omega$ which does not depend on $t$
and $\frac{\partial u}{\partial \nu}$ denotes the directional derivative of $u$
along the outwards pointing unit normal of $\Omega$.
We remark that this setting includes Neumann boundary conditions for $\beta \equiv 0$.
A \emph{mild solution} of~\eqref{diffeq} is a function
$u \in \mathrm{C}\bigl( [0,\infty) ; \mathrm{C}(\overline{\Omega}) \bigr)$
such that $\int_0^t u(s) \dx s \in D(\rlaplace)$ and
\[
	u(t) = u_0 + \rlaplace \int_0^t u(s) \dx s
\]
for every $t \ge 0$.
Note that we have incorporated the Robin boundary conditions
\begin{equation}\label{robinbc}
	\frac{\partial u}{\partial\nu} + \beta u = 0 \text{ on } \partial\Omega
\end{equation}
into the domain
\begin{align*}
	D(\rlaplace) & := \left\{ u \in H^1(\Omega) \cap \mathrm{C}(\overline{\Omega}) \mid \laplace u \in \mathrm{C}(\overline{\Omega}), \vphantom{\int_\Omega} \right. \\
			& \qquad\quad \left. \int_\Omega \nabla u \nabla \phi + \int_\Omega \laplace u \, \phi + \int_{\partial\Omega} u \, \phi \; \beta \dx \sigma = 0 \text{ for every } \phi \in H^1(\Omega) \right\}, \\
	\rlaplace u & := \laplace u
\end{align*}
of the \emph{Laplacian on $\Omega$ subject to Robin boundary conditions}.

It is known that for every non-negative, bounded, measurable function
$\beta$ and every initial value $u_0 \in \mathrm{C}(\overline{\Omega})$
there exists a unique mild solution to problem~\eqref{diffeq}.
In fact, Warma proved that under the above assumptions
$\rlaplace$ generates a $\mathrm{C}_0$-semigroup
$T_R(t)$ on $\mathrm{C}(\overline{\Omega})$~\cite{Warma06},
and it follows from the general theory of semigroups that
then $u(t) = T_R(t) u_0$ is the unique mild solution of the
corresponding homogeneous abstract Cauchy problem~\cite[Proposition~II.6.4]{Nagel00};
note that Warma's proof remains valid for arbitrary
non-negative functions $\beta \in L^\infty(\partial\Omega)$.

If we want to calculate this solution numerically,
a typical problem is how to handle the boundary conditions.
To fix the general ideas, let $\Omega = (0,1)$ and assume that we
want to apply an explicit finite difference method. Then one replaces the
derivatives $u_t$ and $u_{xx}$ by appropriate difference quotients
and successively calculates approximations $u(t_n,x_j)$ of the exact
solution $u(t,x)$ by the relation
\[
	\frac{u(t_{n+1},x_j) - u(t_n,x_j)}{k} = \frac{u(t_n,x_{j+1}) - 2u(t_n,x_j) + u(t_n,x_{j-1})}{h^2},
\]
where $t_n = n \cdot k$ and $x_j = j \cdot h$ for given small numbers $k, h > 0$.
Note, however, that this cannot be directly applied to calculate $u(t_{n+1}, 0)$
and $u(t_{n+1}, 1)$ because $u(t_n, -h)$ and $u(t_n, 1+h)$ are not defined.
For Dirichlet boundary conditions, we can assign $u(t_n, -h) \coloneqq u(t_n, 1+h) \coloneqq 0$.
On the other hand, for Neumann boundary conditions
the situation is not as simple. One common technique is to use
\begin{equation}\label{reflection}
	u(t_n, -h) \coloneqq u(t_n, h) \quad \text{and} \quad u(t_n, 1+h) \coloneqq u(t_n, 1-h)
\end{equation}
in the calculations, which comes from a second order accurate approximation
of the derivative at the boundary, see~\cite[Section~8.3]{Str04}.

The aim of the article at hand is to justify the use of~\eqref{reflection}
for Neumann boundary conditions from a semigroup perspective, showing that the
corresponding continuous method converges to the exact solution as
$k \to 0$, and to extend it to the more general case of Robin boundary 
conditions.
More precisely, we construct an extension operator $E_\beta$ from
$\mathrm{C}(\overline{\Omega})$ to $\mathrm{C}_0(\mathds{R}^N)$ which
depends only on $\Omega$ and $\beta$ (but not on $t$) and resembles
a continuous version of~\eqref{reflection} if $\beta = 0$,
such that $E_\beta$ is a contraction and that $E_\beta u$ is sufficiently
regular whenever $u \in D(\rlaplace)$; we refer to Corollary~\ref{DCC}
for the precise statement.
For operators $E_\beta$ satisfying these two assumptions, we prove the
Trotter-like (compare to~\cite{Trotter59}) approximation result
\begin{equation}\label{claim}
	T_R(t)u = \lim_{n \to \infty} \bigl( R G_0({\textstyle \frac{t}{n}}) E_\beta \bigr)^n u
	\text{ in } \mathrm{C}(\overline{\Omega})
	\text{ for every } u \in \mathrm{C}(\overline{\Omega})
\end{equation}
uniformly on $[0,T]$ for every $T > 0$,
where $G_0(t)$ denotes the Gaussian semigroup on $\mathrm{C}_0(\mathds{R}^N)$
and $R\colon \mathrm{C}_0(\mathds{R}^N) \to \mathrm{C}(\overline{\Omega})$ is the
\emph{restriction operator} $Ru = u|_\Omega$.
This shows how
Robin boundary conditions can be incorporated into a numerical solver such
that the numerical solutions converge uniformly on $[0,T] \times \overline{\Omega}$,
at least if error introduced by space discretization is neglected.

The article is organized as follows.
In Section~\ref{PreRes} we show how to represent a neighborhood of
$\partial\Omega$ in terms of the outwards pointing unit normal and recall some facts
about the Laplacian.
In Section~\ref{ExtOp} we construct the extension operator $E_\beta$
related to the Robin Laplacian
and prove its aforementioned two properties which ensure~\eqref{claim} as we show
in Section~\ref{Approximation}.
Section~\ref{DirichletSec} deals with the limiting case $\beta \to \infty$
giving rise to Dirichlet boundary conditions.
Finally, Section~\ref{OpenProblems} summarizes the results.

\section{Notation and Preliminary Results}\label{PreRes}
It is well-known that for smooth boundary a neighborhood of $\partial\Omega$ can be
para\-metrized by the outwards pointing unit normal $\nu$.
Because certain features of the parametrization are needed later on,
we state this result in the formulation we want to use and prove it
for the sake of completeness.
\begin{proposition}\label{parametrization}
	Let $\delta > 0$ and
	\[
		T: \partial\Omega \times (-\delta, \delta) \to \mathds{R}^N, \; (p, t) \mapsto p + t \nu(p).
	\]
	For small $\delta > 0$, $T$ is a smooth diffeomorphism onto a
	neighborhood of $\partial\Omega$ in $\mathds{R}^N$.
\end{proposition}

\begin{remark}\mbox{}\vspace{-0.5em}
	\begin{enumerate}[(a)]
	\item
		Let $x_0 \in \partial\Omega$ be arbitrary.
		By definition of ``smooth boundary'', $\Omega$ can locally at $x_0$
		be represented as the subgraph of a smooth function
		$\phi\colon U \to \mathds{R}$ ($U \subset \mathds{R}^{N-1}$) up to rotation.
		Assume for the moment that no rotation is needed. Then
		\begin{equation}\label{localboundary}
			\partial \Omega \cap V = \left\{ \begin{pmatrix} z \\ \phi(z) \end{pmatrix} \mid z \in U \right\}
		\end{equation}
		for an open set $V \subset \mathds{R}^N$.
		Thus $z \mapsto \left(\begin{smallmatrix} z \\ \phi(z) \end{smallmatrix}\right)$
		is a bijection of an open subset of $\mathds{R}^{N-1}$ onto a neighborhood of
		$x_0$ in $\partial\Omega$.
		Using these mappings for all $x \in \partial\Omega$ as charts we make $\partial\Omega$
		into a smooth manifold built upon the subspace topology induced by $\mathds{R}^N$.
		Thus we can look at $T$ as a mapping from a manifold to a euclidean space.
		As usual we say that $T$ is smooth if the composition $T^\ast$ of $T$ with a chart is smooth, i.e., if
		\begin{equation}\label{Tast}
			T^\ast: U \times (-\delta, \delta) \to \mathds{R}^N, \;
				(z,t) \mapsto \begin{pmatrix} z \\ \phi(z) \end{pmatrix} + t \nu\!\begin{pmatrix} z \\ \phi(z) \end{pmatrix}
		\end{equation}
		is smooth as a mapping between euclidean spaces.
	\item
		Using charts, the outwards pointing unit normal $\nu$ can be written as
		\begin{equation}\label{normal}
			\nu\!\begin{pmatrix} z \\ \phi(z) \end{pmatrix}
				= \pm\left| \begin{pmatrix} -\nabla \phi(z) & 1 \end{pmatrix}^T\right|^{-1} \begin{pmatrix} -\nabla \phi(z) & 1 \end{pmatrix}^T.
		\end{equation}
		To see this, note that for $x \in \partial\Omega$ the direction of $\nu(x)$ is
		uniquely described by the property
		that for every smooth curve $\xi$ in $\partial\Omega$ satisfying
		$\xi(0) = x = \left(\begin{smallmatrix} z \\ \phi(z) \end{smallmatrix}\right)$ the vectors
		$\nu(x)$ and $\xi'(0)$ are orthogonal. Since by~\eqref{localboundary}, locally
		$\xi(t) = \left(\begin{smallmatrix} \zeta(t) \\ \phi(\zeta(t)) \end{smallmatrix}\right)$, where
		$\zeta(0) = z$, \eqref{normal} follows from the identity
		\[
			\begin{pmatrix} -\nabla \phi(z) & 1 \end{pmatrix} \begin{pmatrix} \zeta'(0) \\ \nabla\phi(\zeta(0)) \cdot \zeta'(0) \end{pmatrix}
				= -\nabla\phi(z) \cdot \zeta'(0) + \nabla\phi(z) \cdot \zeta'(0) = 0.
		\]
	\end{enumerate}
\end{remark}

\begin{proof}[Proof of Proposition~\ref{parametrization}]
	Let $x \in \partial\Omega$ be arbitrary. Working locally near $x$, for simplicity we may assume without
	loss of generality that there exists $\phi$ be as in the previous remark, i.e.,
	that no rotation is needed for $\Omega$ to be the subgraph of a smooth function.
	Then $x = \left(\begin{smallmatrix} z \\ \phi(z) \end{smallmatrix}\right)$
	for some $z$. Using~\eqref{normal}, the derivative of $T^\ast$ is
	\begin{equation}\label{Tastderiv}
		{T^\ast}'(z, 0) = \begin{pmatrix} \Id & -c\nabla\phi(z)^T \\ \nabla\phi(z) & c \end{pmatrix}.
	\end{equation}
	Here, $c = \pm\left| \left(\begin{smallmatrix} -\nabla \phi(z) \\ 1 \end{smallmatrix}\right)\right|^{-1} \neq 0$.
	In particular, we obtain
	\[
		\det {T^\ast}'(z,0)
			= c \cdot \det\begin{pmatrix} \Id & -\nabla\phi(z)^T \\ \nabla\phi(z) & 1 \end{pmatrix}
			= c \cdot \det\begin{pmatrix} \Id & -\nabla\phi(z)^T \\ 0 & 1 + \left| \nabla\phi(z) \right|^2 \end{pmatrix}
			\neq 0
	\]
	by applying the Gauss-Jordan elimination algorithm.
	The inverse function theorem asserts that $T^\ast$ and hence $T$ is locally a smooth diffeomorphism.
	Because $x \in \partial\Omega$ was arbitrary,
	all that remains to show is that $T$ is injective if $\delta > 0$ is small enough.

	By the above argument for every $x \in \partial\Omega$ there exists an open
	neighborhood $O_x$ of $x$ in $\partial\Omega$ and $\delta_x > 0$ such that
	$T$ is a smooth diffeomorphism from $O_x \times (-\delta_z, \delta_z)$ to a
	neighborhood of $x$ in $\mathds{R}^N$.
	By compactness of $\partial\Omega$ we can choose finitely many such $O_{x_i}$ ($i=1,\dots,m$)
	which already cover $\partial\Omega$.
	It is easily proved by contradiction that we can find
	$\delta > 0$ such that for every $x \in \partial\Omega$ there exists an index
	$k(x) \in \{1, \dots, m\}$ with the property that
	$B_{4\delta}(x) \cap \partial\Omega \subset O_{x_{k(x)}}$,
	where $B_r(a)$ denotes the open ball with center $a$ and radius $r$.
	We pick $\delta$ such that $\delta < \delta_{x_i}$ for all $i=1,\dots,m$.

	For this choice of $\delta$, $T$ is injective.
	To see this, let $T(y_1, t_1) = T(y_2, t_2)$ where $y_1, y_2 \in \partial\Omega$
	and $t_1, t_2 \in (-\delta, \delta)$. We estimate
	\[
		0 = \left|T(y_1, t_1) - T(y_2, t_2)\right| \ge \left|y_1 - y_2\right| - \left( \left|t_1\nu(y_1)\right| + \left|t_2\nu(y_2)\right| \right) \ge \left|y_1 - y_2\right| - 2\delta.
	\]
	This shows $|y_1 - y_2| \le 2\delta$ and thus $y_2 \in B_{4\delta}(y_1)$,
	hence $y_1, y_2 \in O_{x_k}$, $k = k(y_2)$. By construction, $T$ is injective
	on $O_k \times (-\delta, \delta)$, hence $y_1 = y_2$ and $t_1 = t_2$,
	proving the claim.
\end{proof}

\begin{lemma}\label{density}
	The set $D := D(\rlaplace) \cap \mathrm{C}^\infty(\overline{\Omega})$
	is an operator core for $\rlaplace$, i.e., $D$ is dense in
	$D(\rlaplace)$ with respect to the graph norm.
\end{lemma}

\begin{proof}
	As $\rlaplace$ is a generator,
	the space $\bigcap_{n \in \mathds{N}} D(\rlaplace^n)$ is a core
	for $\rlaplace$~\cite[Proposition~II.1.8]{Nagel00}. Moreover,
	\[
		R(1, \rlaplace) \left( H^m(\Omega) \cap \mathrm{C}(\overline{\Omega}) \right)
			\subset H^{m+2}(\Omega)
	\]
	for every $m \in \mathds{N}_0$ by the regularization properties of
	elliptic operators~\cite[Remark~2.5.1.2]{Grisvard85}.
	By a standard Sobolev embedding theorem~\cite[Section~5.6]{Evans98},
	\[
		D(\rlaplace^n) = R(1, \rlaplace)^n \mathrm{C}(\overline{\Omega})
			\subset H^{2n}(\Omega)
			\subset \mathrm{C}^{2n - \left[\frac{N}{2}\right] - 1}(\overline{\Omega})
	\]
	for all $n > \frac{N}{4}$. Letting $n \to \infty$ we obtain the assertion.
\end{proof}
We remark that for $u \in \mathrm{C}^\infty(\overline{\Omega})$ the normal derivative
exists in the classical sense. For these functions, $u \in D(\rlaplace)$ is equivalent
to~\eqref{robinbc}, and we will frequently use the boundary condition in this way.

Let $G_2(t)$ denote the $\mathrm{C}_0$-semigroup on $L^2(\mathds{R}^N)$ with generator
\begin{align*}
	D(\Llaplace) & := \left\{ u \in L^2(\mathds{R}^N) \mid \laplace u \in L^2(\mathds{R}^N) \right\}, \\
	\Llaplace u & := \laplace u.
\end{align*}
The semigroup $G_2(t)$ leaves the space $\mathrm{C}_0(\mathds{R}^N) \cap L^2(\mathds{R}^N)$
invariant, and its restriction extends continuously to a positive, contractive $\mathrm{C}_0$-semigroup
on $\mathrm{C}_0(\mathds{R}^N)$, denoted by $G_0(t)$. The generator of this semigroup is
\begin{align*}
	D(\Claplace) & := \left\{ u \in \mathrm{C}_0(\mathds{R}^N) \mid \laplace u \in \mathrm{C}_0(\mathds{R}^N) \right\}, \\
	\Claplace u & := \laplace u.
\end{align*}
We will refer to both semigroups as \emph{the Gaussian semigroup}.
For more details about the Gaussian semigroup, we refer to~\cite[Chapter~3.7]{ABHN01}.

\section{Extension Operator}\label{ExtOp}
Given a smooth bounded open set $\Omega$ and a smooth function
$\beta\colon \partial\Omega \to \mathds{R}_+$, we construct an extension operator $E_\beta$
which satisfies the assumptions under which we will prove~\eqref{claim} in
Section~\ref{Approximation}. For $\beta = 0$, the operator is similar to, but slightly simpler
than the extension operator in~\cite[Section~II.5.4]{Evans98}. However,
the properties which we prove here may also be of independent interest.

For the whole section, let $\delta$ and $T$ be as in Proposition~\ref{parametrization}.
We start by fixing a ``kinking function'' $\rho$. First choose a function
$\rho_1$ having the following properties.
\begin{enumerate}[(a)]
\item
	$\rho_1 \in \mathrm{C}^\infty( [0, \infty) \times [0, \infty) )$
\item\label{contractiveness}
	$0 \le \rho_1(\gamma, t) \le 1$ for all $\gamma, t \ge 0$
\item\label{smallsupport}
	$\rho_1(\gamma, t) = 0$ for all $t \ge \frac{\delta}{2}$ and $\gamma \ge 0$
\item\label{normed}
	$\rho_1(\gamma, 0) = 1$ for all $\gamma \ge 0$
\item
	$\frac{\partial}{\partial t} \rho_1(\gamma, 0) = -2\gamma$ for all $\gamma \ge 0$
\item
	$\frac{\partial^2}{\partial t^2} \rho_1(\gamma, 0) = 4\gamma^2$ for all $\gamma \ge 0$
\end{enumerate}
Here $\frac{\partial}{\partial t} \rho_1$ denotes the partial derivative of
$\rho_1$ with respect to the second argument.
For example, we may choose
$\rho_1(\gamma, t) \coloneqq \exp\left(-2\gamma t\right) \chi(t)$
where $\chi$ is a smooth cut-off function such that
$\chi \equiv 1$ near $0$.

Now define $\rho\colon \Omega^C \to \mathds{R}$ to be
\[
	\rho(x) := \begin{cases}
		\rho_1(\beta(z), t), & \mbox{if } x = T(z, t), \, 0 \le t < \delta, \\
		0, & \mbox{otherwise}. \end{cases}
\]
Note that $\rho$ is well-defined since $T$ is injective,
and it is smooth by construction.

\begin{definition}[Reflection at the Boundary]
	Let $x \in T(\partial\Omega \times (-\delta, \delta))$, $x = T(z, t)$.
	We call $\tilde{x} \coloneqq T(z, -t)$ the
	\emph{(orthogonal) reflection of $x$ at the boundary $\partial \Omega$}.
	For a function $u\colon \Omega \to \mathds{R}$ we define the
	\emph{reflected function}
	\[
		\tilde{u}: T(\partial\Omega \times (0, \delta)) \to \mathds{R}, \quad \tilde{u}(x) \coloneqq u(\tilde{x}).
	\]
	We define the \emph{extension operator $E_\beta$ belonging to $\beta$} as
	\begin{equation}\label{Ebeta}
		E_\beta\colon \mathrm{C}(\overline{\Omega}) \to \mathrm{C}_0(\mathds{R}^N), \quad
		E_\beta u \coloneqq \begin{cases}
			u, & \mbox{on } \overline{\Omega}, \\
			\rho \tilde{u}, & \mbox{on } \overline{\Omega}^C.
		\end{cases}
	\end{equation}
	Here $\rho \tilde{u}$ is understood to be $0$ outside $T(\partial\Omega \times (0,\delta))$
	because $\rho$ equals $0$ in that region.
\end{definition}

\begin{lemma}\label{EbetaProp}
	The operator $E_\beta$ is well-defined, linear, positive, contractive and
	an extension operator, i.e., $RE_\beta = \Id$, where
	$R\colon \mathrm{C}_0(\mathds{R}^N) \to \mathrm{C}(\overline{\Omega})$, $u \mapsto u|_\Omega$.
\end{lemma}

\begin{proof}
	Let $u \in \mathrm{C}(\overline{\Omega})$. 
	By property~\eqref{normed}, the function $E_\beta u$ is continuous on $\mathds{R}^N$.
	Since it has compact support, $E_\beta u \in \mathrm{C}_0(\mathds{R}^N)$.
	Positivity and contractivity follow from property~\eqref{contractiveness}.
	The other two properties are obvious from~\eqref{Ebeta}.
\end{proof}

We now turn towards a more interesting property of $E_\beta$: we prove that
it maps $D$ as defined Lemma~\ref{density} into $D(\Claplace)$.
This extensive calculation is split into several lemmata.
Most calculations will be carried out in local
coordinates, i.e., locally at $x_0 = T^\ast(z_0,0) \in \partial\Omega$,
where we represent all functions with respect to the charts as follows.
Here $T^\ast$ is defined as in~\eqref{Tast}.
\begin{align*}
	u^\ast(z, t) & := u(T^\ast(z, t)) &
	\tilde{u}^\ast(z, t) & := \tilde{u}(T^\ast(z, t)) \\
	\beta^\ast(z) & := \beta(T^\ast(z, 0)) &
	\rho^\ast(z, t) & := \rho(T^\ast(z, t)) = \rho_1(\beta^\ast(z), t)
\end{align*}

In the following we will adhere to the usual notation for normal
derivatives, i.e., $\frac{\partial g}{\partial \nu}$ denotes the
directional derivative of $g$ along the outwards pointing unit normal
\emph{with respect to the domain of $g$}.
Note that for functions defined on $\Omega^C$ this means that
$\frac{\partial g}{\partial \nu} = -\nabla g \cdot \nu$, where $\nu$
always denotes the outwards pointing unit normal of $\Omega$.

\begin{lemma}\label{DL2}
	Let $u \in D$. Then $E_\beta u \in D(\Llaplace)$ and $\bigl(\Llaplace(Eu)\bigr)|_\Omega = \rlaplace u$.
\end{lemma}

\begin{proof}
	The continuous function $E_\beta u$ has compact support, hence $E_\beta u \in L^2(\mathds{R}^N)$.
	Moreover, $E_\beta u$ is smooth away from $\partial\Omega$ being the composition of smooth
	functions. Thus the measurable function
	\[
		f \coloneqq \begin{cases}
				\laplace u, & \mbox{on } \Omega, \\
				\laplace (\rho\tilde{u}), & \mbox{on } \overline{\Omega}^C.
			\end{cases}
	\]
	is defined outside $\partial\Omega$ which is a set of measure zero.
	As $u$ and $\rho \tilde{u}$ are smooth up to $\partial\Omega$,
	$f$ is bounded.
	Note that $f$ has compact support, hence $f \in L^2(\mathds{R}^N)$.

	For the assertion of the lemma, it remains to show
	that $f = \laplace (E_\beta u)$ in the sense of distributions.
	For this we calculate the (classical) normal derivative
	of $\rho \tilde{u}$ using that $u$ satisfies~\eqref{robinbc}.
	For $z \in \partial\Omega$ we have
	\[
		\frac{\partial \tilde{u}}{\partial \nu}(z)
			= -\lim_{h \to 0} \frac{ \tilde{u}(z + h\nu(z)) - \tilde{u}(z) }{h}
			= -\lim_{h \to 0} \frac{u(z - h\nu(z)) - u(z)}{h}
			= -\beta(z) u(z)
	\]
	and
	\[
		\frac{\partial \rho}{\partial \nu}(z)
			= -\lim_{h \to 0} \frac{\rho(z + h\nu(z)) - \rho(z)}{h}
			= -\lim_{h \to 0} \frac{\rho_1(\beta(z), h) - \rho_1(\beta(z), 0)}{h}
			= 2\beta(z).
	\]
	This implies
	\[
		\frac{\partial (\rho\tilde{u})}{\partial \nu}(z)
			= \rho(z) \, \frac{\partial \tilde{u}}{\partial \nu}(z) + \frac{\partial \rho}{\partial \nu}(z) \, \tilde{u}(z)
			= \frac{\partial u}{\partial \nu}(z) + \frac{\partial \rho}{\partial \nu}(z) \, u(z)
			= \beta(z) u(z).
	\]
	Now let $\phi \in \mathcal{D}(\mathds{R}^N)$ be an arbitrary test function.
	From the above calculations, the classical Green formula~\cite[Section~II.1.3]{DL90}
	and $\rho \tilde{u}|_{\partial\Omega} = u|_{\partial\Omega}$, we obtain
	\begin{align*}
		\int_{\mathds{R}^N} (E_\beta u) \laplace \phi
			& = \int_\Omega u \, \laplace \phi + \int_{\Omega^C} \rho \tilde{u} \, \laplace \phi \\
			& = \int_\Omega \laplace u \phi + \int_{\left(\partial \Omega\right)^+} \left(u \frac{\partial \phi}{\partial \nu} - \frac{\partial u}{\partial \nu} \phi \right) \dx\sigma \\
				& \qquad + \int_{\Omega^C} \laplace (\rho\tilde{u}) \phi + \int_{\left(\partial \Omega\right)^-} \left(\rho\tilde{u} \frac{\partial \phi}{\partial \nu} - \frac{\partial (\rho\tilde{u})}{\partial \nu} \phi \right) \dx\sigma \\
			& = \int_\Omega \laplace u \, \phi + \int_{\Omega^C} \laplace(\rho \tilde{u}) \, \phi
			= \int_{\mathds{R}^N} f \, \phi,
	\end{align*}
	where $(\partial\Omega)^+$ is understood as the (oriented) boundary of $\Omega$, whereas
	$(\partial\Omega)^-$ denotes the (oriented) boundary of $\smash{\overline{\Omega}}^C$.
	This shows $\laplace(E_\beta u) = f$ in the sense of distributions.
\end{proof}

\begin{remark}\label{whattoshow}
	The above lemma tells us that $\laplace(E_\beta u)$ is a function.
	To see that $E_\beta u \in D(\Claplace)$, it remains to show that
	$\laplace(E_\beta u) \in \mathrm{C}_0(\mathds{R}^N)$.
	We already know that $\laplace(E_\beta u)$ has compact support and is smooth on
	$\mathds{R}^N \setminus \partial\Omega$.
	Thus it suffices to show that the function
	can continuously be extended to $\partial\Omega$.
	This is a local property.
	In fact, since we already know that the limits from the inside and the
	outside both exist, it suffices to show that
	$\laplace u(x_0) = \laplace (\rho \tilde{u})(x_0)$
	for every $x_0 \in \partial\Omega$.
\end{remark}

Let $x_0 \in \partial\Omega$ be fixed.
To simplify notation, we may  assume that $\nu(x_0) = e_N$
without loss of generality,
exploiting the rotational invariance of the Laplacian.
Here and in the following, $e_n$ denotes the
the $n$\textsuperscript{th} unit vector in $\mathds{R}^N$.
Moreover, since we treat the problem locally, we may
work in local coordinates, $x_0 = T^\ast(z_0, 0)$.
We start by calculating the partial derivatives of ${T^\ast}^{-1}$,
where $T^\ast$ is defined as in~\eqref{Tast}.
\begin{lemma}\label{Tinvderiv}
	For $n \in \{ 1, \dots, N \}$,
	\begin{align*}
		\left(\frac{\partial}{\partial x_n} \left( {T^\ast}^{-1} \right) \right)(x_0) & = e_n, \\
		\left(\frac{\partial^2}{\partial x_n^2} \left( {T^\ast}^{-1} \right) \right)(x_0) & =
			\begin{cases} \begin{pmatrix} 0 \\ -\frac{\partial^2}{\partial z_n^2} \phi(z_0) \end{pmatrix} & \mbox{if } n \neq N, \\
			0 & \mbox{if } n = N. \end{cases}
	\end{align*}
\end{lemma}

\begin{proof}
	The assumption $\nu(z_0) = e_N$ implies $\nabla\phi(z_0) = 0$ due to~\eqref{normal}.
	As in~\eqref{Tastderiv}, this shows ${T^\ast}'(z_0, 0) = \Id$. By the inverse function theorem,
	\[
		\left( {T^\ast}^{-1} \right)'(x) = \left( {T^\ast}'\left({T^\ast}^{-1}(x)\right) \right)^{-1}.
	\]
	For the partial derivatives at $x_0$ this means
	\[
		\left(\frac{\partial}{\partial x_n} \left( {T^\ast}^{-1} \right) \right)(x_0)
			= {T^\ast}'(z_0, 0)^{-1} e_n = \Id e_n = e_n.
	\]

	To calculate the second derivatives, we employ a differentiation rule
	for matrices,
	$\frac{\mathrm{d}}{\mathrm{d}t} \left( A(t)^{-1} \right) = -A^{-1}(t) A'(t) A^{-1}(t)$.
	\begin{align*}
		& \left(\frac{\partial^2}{\partial x_n^2} \left( {T^\ast}^{-1} \right) \right)(x)
			= \left( \frac{\partial}{\partial x_n} \left( {T^\ast}'\left({T^\ast}^{-1}(x)\right)\right)^{-1} \right) e_n \\
			& = -\left( {T^\ast}'\left({T^\ast}^{-1}(x)\right)\right)^{-1}
					\left( \frac{\partial}{\partial x_n} \left( {T^\ast}'\left({T^\ast}^{-1}(x)\right)\right)\right)
					\left( {T^\ast}'\left({T^\ast}^{-1}(x)\right)\right)^{-1} e_n
	\end{align*}
	If we denote the entries of ${T^\ast}'$ by $t_{ij}$ ($i,j=1,\dots,N$),
	we can proceed as follows.
	\[
		\frac{\partial}{\partial x_n} \left(t_{ij}\left({T^\ast}^{-1}(x)\right)\right)
			= \nabla t_{ij}\left({T^\ast}^{-1}(x)\right) \cdot
				\left(\frac{\partial}{\partial x_n} \left( {T^\ast}^{-1} \right) \right)(x)
	\]
	For $x = x_0$ this yields
	\[
		\frac{\partial}{\partial x_n} t_{ij}\left({T^\ast}^{-1}(x_0)\right)
			= \nabla t_{ij}(z_0, 0) e_n
			= \frac{\partial}{\partial z_n} t_{ij}(z_0, 0),
	\]
	where for notational simplicity we use $z_N$ as an alias for the variable $t$.
	Inserting this expression into the above identity, we arrive at
	\begin{align*}
		\left(\frac{\partial^2}{\partial x_n^2} \left( {T^\ast}^{-1} \right) \right)(x_0)
			& = -\left( \frac{\partial}{\partial z_n} t_{ij}(z_0, 0) \right)_{i,j=1,\dots,N} e_n \\
			& = -\left( \frac{\partial}{\partial z_n} t_{in}(z_0, 0) \right)_{i=1,\dots,N}
			= -\frac{\partial^2}{\partial z_n^2} T^\ast(z_0, 0).
	\end{align*}
	In combination with formula~\eqref{Tast} this finishes the proof.
\end{proof}

Having the derivatives of the charts at hand, we are able to calculate
all derivatives in local coordinates.
\begin{lemma}
	Let $f$ and $f^\ast$ be functions such that locally
	$f^\ast(z, t) = f(T(z, t))$. Then
	\begin{align*}
		\nabla f(x_0) & = \nabla f^\ast(z_0, 0), \\
		\laplace f(x_0)
			& = \sum_{n=1}^{N-1} \frac{\partial^2}{\partial z_n^2} f^\ast(z_0, 0)
			+ \frac{\partial^2}{\partial t^2} f^\ast(z_0, 0)
			- \frac{\partial}{\partial t} f^\ast(z_0, 0) \sum_{n=1}^{N-1} \frac{\partial^2}{\partial z_n^2} \phi(z_0).
	\end{align*}
	In particular,
	\begin{align*}
		\nabla \rho(x_0) & = \begin{pmatrix} 0 & -2\beta(x_0) \end{pmatrix}, &
		\laplace \rho(x_0) & = 4\beta(x_0)^2 + 2\beta(x_0) \sum_{n=1}^{N-1} \frac{\partial^2}{\partial z_n^2} \phi(z_0).
	\end{align*}
\end{lemma}

\begin{proof}
	Differentiating $f(x) = f^\ast\left({T^\ast}^{-1}(x)\right)$ we obtain
	\begin{align*}
		\frac{\partial}{\partial x_n} f(x)
			& = \left( \nabla f^\ast \right) \left( {T^\ast}^{-1}(x) \right) \left(\frac{\partial}{\partial x_n} \left( {T^\ast}^{-1} \right) \right)(x), \\
		\frac{\partial^2}{\partial x_n^2} f(x)
			& = \left(\frac{\partial}{\partial x_n} \left( {T^\ast}^{-1} \right) \right)^T\!\!(x) \, H_{f^\ast}\left( {T^\ast}^{-1}(x) \right)
				\left(\frac{\partial}{\partial x_n} \left( {T^\ast}^{-1} \right) \right)(x)  \\
			& \quad + \left( \nabla f^\ast \right) \left( {T^\ast}^{-1}(x) \right) \left(\frac{\partial^2}{\partial x_n^2} \left( {T^\ast}^{-1} \right) \right)(x),
	\end{align*}
	where $H_{f^\ast} = \left( \frac{\partial^2}{\partial z_i \, \partial z_j} f^\ast \right)_{i,j=1,\dots,N}$
	denotes the Hessian matrix of $f^\ast$.
	By using Lemma~\ref{Tinvderiv} and summing up, we arrive at the desired formulae for $x = x_0$.

	Concerning $\rho$ we remark that $\rho^\ast(z, 0) = \rho(\beta(z), 0) = 1$
	implies $\frac{\partial}{\partial z_n} \rho^\ast(z_0, 0) = 0$ ($n=1,\dots,N-1$).
	On the other hand, the derivatives with respect to $t$ equal
	\begin{align*}
		\frac{\partial}{\partial t} \rho^\ast(z_0, 0) & = \frac{\partial}{\partial t} \rho_1(\beta^\ast(z_0), 0) = -2\beta^\ast(z_0) = -2\beta(x_0), \\
		\frac{\partial}{\partial t^2} \rho^\ast(z_0, 0) & = \frac{\partial}{\partial t^2} \rho_1(\beta^\ast(z_0), 0) = 4\beta^\ast(z_0)^2 = 4\beta(x_0)^2.
	\end{align*}
	With this information, the formulae for $\rho$ follow from the general formulae.
\end{proof}

Finally, it is easy to calculate the relation between the derivatives of the
function and its reflection at the boundary in local coordinates.
It suffices to observe that
\[
	\tilde{u}^\ast(z, t) = \tilde{u}(T(z, t)) = u(T(z, -t)) = u^\ast(z, -t).
\]
From this we deduce the following formulae.
\begin{align*}
	\tilde{u}^\ast(z, t) & = u^\ast(z, -t) \\
	\frac{\partial}{\partial z_n} \tilde{u}^\ast(z, t) & = \frac{\partial}{\partial z_n}u^\ast(z, -t) &
	\frac{\partial}{\partial t} \tilde{u}^\ast(z, t) & = -\frac{\partial}{\partial t}u^\ast(z, -t) \\
	\frac{\partial^2}{\partial z_n^2} \tilde{u}^\ast(z, t) & = \frac{\partial^2}{\partial z_n^2}u^\ast(z, -t) &
	\frac{\partial^2}{\partial t^2} \tilde{u}^\ast(z, t) & = \frac{\partial^2}{\partial t^2}u^\ast(z, -t)
\end{align*}

Now we are ready to prove continuity of $\laplace(E_\beta u)$ at $x_0$.
\begin{proposition}\label{contpoint}
	For every $u \in D$, $\laplace u(x_0) = \laplace(\rho \tilde{u})(x_0)$.
\end{proposition}

\begin{proof}
	Note that
	\[
		\frac{\partial}{\partial t} \tilde{u}^\ast(z_0, 0)
			= -\frac{\partial}{\partial t} u^\ast(z_0, 0)
			= -\frac{\partial u}{\partial \nu}(x_0)
			= \beta(x_0) u(x_0)
			= \beta(x_0) \tilde{u}(x_0).
	\]
	We use the formulae of this section to obtain the desired identity.
	\begin{align*}
		\laplace(\rho \tilde{u})(x_0)
			& = \laplace \rho(x_0) \; \tilde{u}(x_0) + 2 \nabla\rho(x_0) \cdot \nabla\tilde{u}(x_0) + \rho(x_0) \; \laplace\tilde{u}(x_0) \\
			& \hspace{-2.5em} = 4\beta(x_0)^2 \tilde{u}(x_0) + 2\beta(x_0) \tilde{u}(x_0) \sum_{n=1}^{N-1} \frac{\partial^2}{\partial z_n^2} \phi(z_0)
				- 4\beta(x_0) \frac{\partial}{\partial t}\tilde{u}^\ast(z_0, 0) \nonumber \\
				& \quad + \sum_{n=1}^{N-1} \frac{\partial^2}{\partial z_n^2} \tilde{u}^\ast(z_0, 0) + \frac{\partial^2}{\partial t^2} \tilde{u}^\ast(z_0, 0)
				- \frac{\partial}{\partial t} \tilde{u}^\ast(z_0, 0) \sum_{n=1}^{N-1} \frac{\partial^2}{\partial z_n^2} \phi(z_0) \\
			& \hspace{-2.5em} = \sum_{n=1}^{N-1} \frac{\partial^2}{\partial z_n^2} u^\ast(z_0, 0) + \frac{\partial^2}{\partial t^2} u^\ast(z_0, t)
				- \frac{\partial}{\partial t} u^\ast(z_0, 0) \sum_{n=1}^{N-1} \frac{\partial^2}{\partial z_n^2} \phi(z_0) \\
			& \hspace{-2.5em} = \laplace u(x_0).
	\end{align*}
\end{proof}

The following theorem is the main result of this section.
As explained in Remark~\ref{whattoshow}, it follows by combining
Lemma~\ref{DL2} and the last proposition.
Even though Theorem~\ref{DC} is also true for the usual
extension operator for Lipschitz domains~\cite[VI.\S 3, Theorem~5]{Stein70},
that operator fails to be contractive and thus
is more difficult to handle for the application
in Section~\ref{Approximation}.
\begin{theorem}\label{DC}
	The operator $E_\beta$ maps $D$ into $D(\Claplace)$.
\end{theorem}
\begin{corollary}\label{DCC}
	The operator $E_\beta$ maps $D(\rlaplace)$ into $D(\Claplace)$.
\end{corollary}
\begin{proof}
	There exists a constant $C > 0$ satisfying
	$\| E_\beta u \|_{D(\Claplace)} \le C \| u \|_{D(\rlaplace)}$
	for all $u \in D$.
	To see this, note that on $\smash{\overline{\Omega}}^C$
	\begin{align*}
		\left\| \laplace(\tilde{u} \rho) \right\|_\infty
			& = \left\| \laplace \tilde{u} \rho + 2 \nabla \tilde{u} \nabla \rho + \tilde{u} \laplace \rho \right\|_\infty \\
			& \le \left\| \laplace \tilde{u} \right\|_\infty \left\| \rho \right\|_\infty + 2 \left\| \nabla \rho \right\|_\infty \left( \eps \left\| \laplace \tilde{u} \right\|_\infty + \left\| \tilde{u} \right\|_\infty \right) + \left\| \tilde{u} \right\|_\infty \left\| \laplace \rho \right\|_\infty
	\end{align*}
	for every $\eps > 0$. Similarly,
	$\left\| \tilde{u} \right\|_\infty \le \left\| u \right\|_\infty$ and
	$\left\| \laplace \tilde{u} \right\|_\infty \le \widetilde{C} \left( \left\| \laplace u \right\|_\infty + \left\| u \right\|_\infty \right)$,
	using the definition of $\tilde{u}$ as a composition of $u$ and a function
	involving $T$, where $\widetilde{C} > 0$ depends only depends on a estimate
	on the derivatives of $T$. Noting that $\rho$ and its derivatives are bounded
	by assumption, we see that there exists a $C > 0$ as claimed.

	As $D$ is a core of $\rlaplace$, the above estimate shows that there exists
	a unique continuous extension of $E_\beta|_D$ to $D(\rlaplace)$, and
	that this operator still takes values in $D(\Claplace)$.
	Because $D(\rlaplace)$ is continuously embedded into $\mathrm{C}(\overline{\Omega})$
	and $E_\beta$ is continuous on $\mathrm{C}(\overline{\Omega})$,
	this extension agrees with $E_\beta|_{D(\rlaplace)}$. Thus the claim
	is proved.
\end{proof}

\section{Approximation Result}\label{Approximation}

In this section, we we prove that
if $E_\beta\colon \mathrm{C}(\overline{\Omega}) \to \mathrm{C}_0(\mathrm{R}^N)$
is a contractive extension operator mapping an operator core $D$
for $\rlaplace$ into $D(\Claplace)$, then
formula~\eqref{claim} holds.
Note that the operator defined in~\eqref{Ebeta} has this properties
as shown in the preceding section.
The tool we use for the proof is the following approximation result
for semigroups due to Chernoff.
\begin{theorem}[{\cite[Theorem~III.5.2]{Nagel00}}]\label{chernoff}
	Let $X$ be a Banach space. Consider a function
	\[
		V: [0, \infty) \to \mathscr{L}(X)
	\]
	satisfying $V(0) = \Id$ and $\|V(t)^m\| \le M$ for all $t \ge 0$, $m \in \mathds{N}$
	and some $M \ge 1$. Assume that
	\[
		Ax \coloneqq \lim_{h \to 0} \frac{V(h)x - x}{h}
	\]
	exists for all on $x \in D \subset X$, where
	$D$ and $(\lambda_0 - A)D$ are dense subspaces in $X$ for some $\lambda_0 > 0$.
	\\
	Then $(A, D)$ is closable and $\overline{A}$ generates a bounded $\mathrm{C}_0$-semigroup $T(t)$,
	which is given by
	\[
		T(t)x = \lim_{n \to \infty} \bigl( V({\textstyle \frac{t}{n}}) \bigr)^n x
	\]
	for every $x \in X$ locally uniformly with respect to $t \ge 0$.
\end{theorem}
We apply the theorem by setting
\begin{align}\label{setting}
	X & \coloneqq \mathrm{C}(\overline{\Omega}), &
	V(t) & \coloneqq RG_0(t)E_\beta, &
	A & \coloneqq \rlaplace.
\end{align}
As $D$ is an operator core for $\rlaplace$,
the density conditions are fulfilled because $\lambda - \rlaplace$
is an isomorphism between $D(\rlaplace)$ with the graph norm and
$\mathrm{C}(\overline{\Omega})$ for every $\lambda > 0$.

\begin{theorem}\label{extcond}
	Let $E_\beta\colon \mathrm{C}(\overline{\Omega}) \to \mathrm{C}_0(\mathds{R}^N)$
	be a contractive extension operator which maps an operator core $D$ for $\rlaplace$
	into $D(\Claplace)$. Then formula~\eqref{claim} holds true.
\end{theorem}

\begin{proof}
	We check the conditions of Chernoff's product formula with the choices made in~\eqref{setting}.
	The fact $V(0) = \Id$ is equivalent to $E_\beta$ being an extension operator.
	Since all three of their factors are contractions, the operators $V(t)$ are contractions for
	every $t \ge 0$, thus $\|V(t)^m\| \le 1$; in particular $V(t)$ is a bounded operator for every $t \ge 0$.
	The density assumptions on $D$ are fulfilled because $D$ is an operator core for $\rlaplace$.

	Now let $u \in D$ be arbitrary. By assumption, $E_\beta u \in D(\Claplace)$. By definition
	of the infinitesimal generator,
	\[
		\frac{V(h)u - u}{h} = R \frac{G(h)(E_\beta u) - E_\beta u}{h} \to R \Claplace E_\beta u. \qquad (h \to 0)
	\]
	Since the function $E_\beta u$ agrees with $u$ on $\Omega$, they represent the same
	distribution acting on the test functions $\mathcal{D}(\Omega)$. This means that they
	have the same distributional derivatives, hence $R \Claplace E_\beta u = \rlaplace u$.
	Having checked all the conditions of Theorem~\ref{chernoff}, we deduce that
	indeed~\eqref{claim} holds true.
\end{proof}

\begin{corollary}
	Formula~\eqref{claim} holds true for the operator $E_\beta$ defined in~\eqref{Ebeta}.
\end{corollary}

\begin{remark}\label{Neumann}
	As a special case, we may choose $\beta = 0$. Then $\rlaplace = \nlaplace$ is the
	Laplacian with \emph{Neumann boundary conditions}. In this case, $E_0$ is the
	reflection without ``kinking'', corresponding to certain numeric schemes
	where Neumann boundary conditions are realized as in~\eqref{reflection}.
	A different extension operator for Neumann boundary conditions
	would be given by extending constantly along the outwards pointing unit normal and again multiplying by a cut-off
	function. This corresponds to a first-order accurate boundary condition approximation,
	see again~\cite[Section~8.3]{Str04}. Although this might seem more natural at first, it is not obvious
	whether formula~\eqref{claim} is true for this choice of $E_\beta$. Unfortunately,
	Chernoff's theorem cannot be applied again because $\laplace E_\beta u$ fails
	to be continuous at $\partial\Omega$ as can easily be seen.
\end{remark}

\section{Dirichlet Boundary Conditions}\label{DirichletSec}
Next we treat the model problem of an elliptic operator on a bounded set,
the Laplacian with Dirichlet boundary conditions. Typically, all results
about elliptic operators are much simpler for this special case.
Surprisingly, for the aim of this article there arise completely different problems
than for Robin and Neumann boundary conditions. This is the reason why we consider it worthwhile
to treat this operator in detail.

The \emph{Laplacian with Dirichlet boundary conditions} defined by
\[
	D(\dlaplace) \coloneqq \left\{ u \in \mathrm{C}_0(\Omega) \mid \laplace u \in \mathrm{C}_0(\Omega) \right\}, \qquad
	\dlaplace u \coloneqq \laplace u
\]
generates an positive, contractive $\mathrm{C}_0$-semigroup $T_D(t)$ on
$\mathrm{C}_0(\Omega)$~\cite[Theorem~6.1.8]{ABHN01}.
Formally, the boundary conditions~\eqref{robinbc} become the \emph{Dirichlet boundary conditions}
$u = 0$ on $\partial\Omega$ in the limit $\beta \to \infty$. This observation can be made
precise, cf.~\cite[Proposition~3.5.3]{Warma02}.
As we want to prove an analogue of~\eqref{claim} for $T_D(t)$, we have
to define an appropriate extension operator $E_\infty$ for $\beta = \infty$.
Taking the limit in~\eqref{Ebeta}, we arrive at
\[
	E_\infty\colon \mathrm{C}_0(\Omega) \to \mathrm{C}_0(\mathds{R}^N), \;
		E_\infty u \coloneqq \begin{cases} u & \text{on } \Omega, \\ 0 & \text{on } \Omega^C. \end{cases}
\]
Note that we had to replace $\mathrm{C}(\overline{\Omega})$ by
$\mathrm{C}_0(\Omega)$ as we require $E_\infty u$ to be continuous.
Unfortunately, we cannot simply replace $E_\beta$ by $E_\infty$ in
formula~\eqref{claim} because the iteration scheme does not
remain in $\mathrm{C}_0(\Omega)$, hence leaving the domain of $E_\infty$.

However, the analogue formula is well-defined (and true) in the $L^2$-context.
To see this, note that $L^2(\Omega)$ is a closed subspace of $L^2(\mathds{R}^N)$
if we consider its functions to be extended by zero. Then the identity mapping
takes the role of $E_\infty$, and the restriction becomes multiplication
with $\setone_\Omega$. Thus, the analogue of formula~\eqref{claim} for Dirichlet
boundary conditions reads
\begin{equation}\label{DirichletApprox}
	T_{D,2}(t) u = \lim_{n \to \infty} \bigl( \setone_\Omega G_2({\textstyle \frac{t}{n}}) \bigr)^n u
	\text{ in } L^2(\Omega) \text{ for every } u \in L^2(\Omega),
\end{equation}
where $T_{D,2}$ denotes the Dirichlet
semigroup on $L^2(\Omega)$ generated by the Laplacian on $L^2(\Omega)$
with domain $H^1_0(\Omega) \cap H^2(\Omega)$.
Indeed, formula~\eqref{DirichletApprox} remains true
even if $\Omega$ has merely Lipschitz regular boundary, cf.~\cite{MS03}.

It is interesting to note that~\eqref{DirichletApprox} cannot be proved using
Chernoff's product formula in the way we did in Section~\ref{Approximation}. For this,
a dense subspace of $H^1_0(\Omega) \cap H^2(\Omega)$ would have to be contained
in $H^2(\mathds{R}^N)$, where both spaces carry the graph norm of the Laplacian.
But then, continuity asserts $H^1_0(\Omega) \cap H^2(\Omega) \subset H^2(\mathds{R}^N)$.
However, this cannot be true. In fact, a function in
$\mathrm{C}_0(\Omega) \cap \mathrm{C}^\infty(\overline{\Omega}) \subset H^1_0(\Omega) \cap H^2(\Omega)$
whose normal derivative does not vanish is not an element of $H^2(\mathds{R}^N)$.

\bigskip
Despite those problems, it is possible to prove a similar result in the same
spirit as in Section~\ref{Approximation} even in $\mathrm{C}_0(\Omega)$.
For this, we need to replace $\setone_\Omega$ by a sequence of smooth
interior cut-off functions. But we have to assure that they exhaust $\Omega$
sufficiently fast compared to the decay of functions in a core for $\dlaplace$.
So we start by a investigation of that decay.
\begin{lemma}\label{boundaryestimate}
	Given a Dirichlet regular bounded set $\Omega$, there exists $m \in \mathds{N}$
	having the following property. Given $t > 0$, there exists a neighborhood
	$U_t$ of $\partial \Omega$ such that the estimate
	\[
		|u(x)| \le t^2 \left\|(\Id - \dlaplace)^m u\right\|_\infty
	\]
	holds for every $x \in \overline{\Omega} \cap U_t$ and every $u \in D(\dlaplace^m)$.
\end{lemma}
\begin{proof}
	It is well-known that $T_D(t)$ has a kernel representation with a
	continuous non-negative symmetric kernel $k_s(x,y)$ which vanishes on $\partial\Omega$ and is
	dominated by the Gaussian kernel~\cite[Section~3.4]{Davies89}.
	Let $m > \frac{N}{2}$ be fixed. The integral formula for powers of the
	resolvent~\cite[Corollary~2.1.11]{Nagel00} shows that $(\Id - \dlaplace)^{-m}$
	is a positive kernel operator with the continuous non-negative symmetric kernel
	\[
		k(x,y) = \int_0^\infty \frac{s^{m-1}}{(m-1)!} \e^{-s} k_s(x,y) \; \dx s
	\]
	which vanishes for $x \in \partial\Omega$. Using compactness of $\overline{\Omega}$
	and $\partial\Omega$ we deduce that for any $\eps > 0$ there exists a
	neighborhood $S_\eps$ of $\partial\Omega$ such that $x \in S_\eps$ implies
	$k(x,y) \le \eps$ for all $y \in \overline{\Omega}$.
	Define $U_t \coloneqq S_\eps$, where $\eps \coloneqq \frac{t^2}{|\Omega|}$.

	Now fix $u \in D(\dlaplace^m)$ and define
	$v \coloneqq (\Id - \dlaplace)^m u \in \mathrm{C}_0(\Omega)$.
	For $x \in \overline{\Omega} \cap U_t$, i.e.,
	$x \in \overline{\Omega} \cap S_\eps$, we obtain
	\[
		|u(x)| \le (\Id - \dlaplace)^{-m} \left|v\right|
			= \int_\Omega k(x,y) |v(y)| \; \dx y
			\le \eps \left|\Omega\right| \left\|v\right\|_\infty
			= t^2 \left\|(\Id - \dlaplace)^m u\right\|_\infty.
	\]
	This concludes the proof.
\end{proof}

We have already explained why we cannot use $E_\infty$ as extension operator.
Instead, we choose
\[
	E_D\colon \mathrm{C}_0(\Omega) \to \mathrm{C}_0(\mathds{R}^N), \;
		E_D u \coloneqq \begin{cases}
			u, & \text{on } \overline{\Omega}, \\
			-\rho \tilde{u}, & \text{on } \smash{\overline{\Omega}}^C,
		\end{cases}
\]
similarly to~\eqref{Ebeta}. Here $\rho$ denotes a cut-off function that equals $1$ near $\partial\Omega$.
Using the same ideas as in Section~\ref{ExtOp} it can be shown that $E_D$
is a contractive extension operator that maps
$D(\dlaplace) \cap \mathrm{C}^\infty(\overline{\Omega})$ into $D(\Claplace)$.
In fact, the main difference to Section~\ref{ExtOp} is that we know
$\laplace u \in \mathrm{C}_0(\Omega)$ for $u \in D(\dlaplace)$ which makes
it easy to check the continuity of $\laplace (E_D u)$, significantly shortening
the chain of arguments.

Now let $m$ be as in the above lemma, and choose a family $(U_t)_{t > 0}$ as in the lemma.
For every $t > 0$ we fix a suitable cut-off function $\chi_t \in \mathrm{C}_0(\Omega)$
satisfying $0 \le \chi_t \le 1$, and $\chi_t(x) = 1$ if $x \in \overline{\Omega} \setminus U_t$.
Moreover, define $\chi_0 \coloneqq \setone_{\overline{\Omega}}$.
To simplify notation, we use the multiplication operator $\chi_t$ as an operator
from $\mathrm{C}_0(\mathds{R}^N)$ to $\mathrm{C}_0(\Omega)$ and $\chi_0$
as the restriction from $\mathrm{C}_0(\mathds{R}^N)$ to $\mathrm{C}(\overline{\Omega})$,
whenever they are applied to functions in $\mathrm{C}_0(\mathds{R}^N)$.

We remark that in view of the kernel of $T_D(t)$ ($t > 0$) being strictly positive
in the interior of $\Omega$ due to the strong maximum principle, it can be seen that
for every compact set $K \subset \Omega$ there exists $t_0 > 0$ such that
$U_t$ and $K$ are disjoint whenever $t < t_0$, implying that
$\chi_t \to \setone_\Omega$ pointwise as $t \to 0$.
In this sense, the next result is another flavor of formula~\eqref{DirichletApprox}.
\begin{theorem}
	Let $m \in \mathds{N}$ and $(\chi_t)_{t \ge 0}$ be as above. Then
	\[
		T_D(t)u = \lim_{n \to \infty} \bigl( \chi_{\frac{t}{n}} \, G_0({\textstyle \frac{t}{n}}) E_D \bigr)^n u
	\]
	for every $u \in \mathrm{C}_0(\Omega)$ uniformly on $[0,T]$ for every $T > 0$.
\end{theorem}
\begin{proof}
	We apply Theorem~\ref{chernoff} to the operators
	\[
		V(t)\colon \mathrm{C}_0(\Omega) \to \mathrm{C}_0(\Omega), \; u \mapsto \chi_t G_0(t) E_D u.
	\]
	The properties $V(0) = \Id$ and $\left\| V(t)^n \right\| \le 1$ for every $t \ge 0$ and
	$n \in \mathds{N}$ are obvious from the properties of $E_D$ and the Gaussian semigroup.
	Let $D \coloneqq D(\dlaplace^m) \cap \mathrm{C}^\infty(\overline{\Omega})$, which is a
	core for $\dlaplace$. This choice makes the density conditions automatic once we show that the
	limit operator is $\dlaplace$.
	
	It only remains to prove the convergence to $\laplace u$ on $D$. For this, let $u \in D$.
	In particular $u \in D(\dlaplace)$, thus $\laplace u \in \mathrm{C}_0(\Omega)$. Note that
	\begin{align*}
		& \left\| \frac{V(t)u - u}{t} - \laplace u \right\|_\infty
			= \left\| \frac{\chi_t G_0(t) E_D u - u}{t} - \laplace u \right\|_\infty \\
			& \quad \le \left\| \chi_t \left( \frac{G_0(t)E_Du - E_Du}{t} - \laplace u \right) \right\|_\infty
				+ \left\| \frac{\chi_t E_Du - u}{t} \right\|_\infty
				+ \left\| \chi_t \laplace u - \laplace u \right\|_\infty.
	\end{align*}
	We estimate the three summands separately. After estimating $\chi_t$ by $1$ in the
	first expression, convergence to zero follows from $E_D u \in D(\Claplace)$ and
	the fact that $\Claplace (E_D u) = \laplace u$ on $\Omega$.
	The third summand can be estimated
	by $\sup_{x \in U_t}  2 \left| \laplace u(x) \right|$ using that $\chi_t = 1$
	on $\overline{\Omega} \setminus U_t$. But since we assumed that $U_t$ leaves any
	compact set $K \subset \Omega$ for small $t$, this expression becomes small as
	$t \to 0$ because $\laplace u \in \mathrm{C}_0(\Omega)$.
	The second summand can be estimated with the help of Lemma~\ref{boundaryestimate}. We obtain
	\[
		\left\| \frac{\chi_t E_D u - u}{t} \right\|_\infty
			= \frac{1}{t} \left\| \chi_t u - u \right\|_\infty
			\le \frac{2}{t} \sup_{x \in U_t} \left| u(x) \right|
			\le 2t \left\| (\Id - \dlaplace)^m u \right\|_\infty
			\to 0
	\]
	as $t \to 0$.
	Together, these three estimates show the convergence of the difference
	quotient to $\laplace u$ as $t$ tends to zero. We have checked the assumptions
	of Chernoff's product formula, thus proving the claim of the theorem.
\end{proof}

\section{Conclusion}\label{OpenProblems}
It is a direct consequence of~\eqref{claim} that $T_R(t)$ is a positive
semigroup. Because the operators on the right are $L^\infty$-contractive, it is
also clear that $T_R(t)$ is $L^\infty$-contractive, thus submarkovian.
In the same way other properties of the limiting semigroup can be deduced
by such an approximation formula, as long as they are preserved when taking
limits in the strong operator topology.
To obtain further properties, it might help to modify the formula a little bit.

So far, we have only considered the Gaussian semigroup as underlying tool. However,
it can be seen from the proofs that actually we used only few properties of
the Gaussian semigroup. More precisely, we only used that $G_0(t)$ is a contraction
on $\mathrm{C}_0(\mathds{R}^N)$ and that any continuous function $u$ such that
the support of $u$ is contained in a given neighborhood of $\Omega$ and
$\laplace u$ is continuous on $\mathds{R}^N$ is in the domain of the generator
of $G_0(t)$. Thus we could replace $G_0(t)$ with other semigroups, for example
with the semigroup generated by the Laplacian with Dirichlet or Neumann
boundary conditions, on a larger bounded set $\Omega' \subset \mathds{R}^N$.
Then~\eqref{claim} becomes an approximation formula where the approximating
operators are compact.
Note, however, that this does not imply that $T_R(t)$ is compact, as the
limit is only in the strong operator topology.

Similarly, we can try to approximate $T_R(t)$ only in terms of operators
on $\mathrm{C}(\overline{\Omega})$, i.e., without any extension to $\mathds{R}^N$,
to obtain an intrinsic approximation. The most natural Trotter-like candidate
of this kind would be
$S_n(t) \coloneqq \bigl( T_D(\frac{\alpha t}{n}) T_N(\frac{\beta t}{n}) \bigr)^n$,
where $\alpha$ and $\beta$ are positive numbers such that $\alpha + \beta = 1$
and $T_D(t)$ and $T_N(t)$ denote the Dirichlet and Neumann semigroups on
$\mathrm{C}(\overline{\Omega})$, respectively. It is known, however, that
$\lim S_n(t) = T_D(t)$ in the strong operator topology on $L^2(\Omega)$,
whenever $\alpha > 0$, see~\cite{Kato78}.

\bigskip

Recall that regarding the extension operator we used only two of its properties
in Section~\ref{Approximation}, namely contractivity and
some regularity of the extended function.
By definition of the extension operator, contractivity
came for free. This is due to the rather special definition
of $E_\beta$ and the choice of spaces and
is a very convenient prerequisite for the application
of Theorem~\ref{chernoff}, although not a necessary one.

Assume that we replace $E_\beta$ by some other, non-contractive
extension operator $E$. This is a natural consideration
because most extension operators are non-contractive.
In fact, it is easy to see that no operator extending
$\mathrm{C}^1(\overline{\Omega})$ to
$\mathrm{C}^1(\mathds{R}^N)$ can be contractive.
This also shows that it is a very special property for
an extension operator to be contractive \emph{and} to preserve
the regularity of functions in $D(\rlaplace)$.

For such an extension operator $E$, it is considerably more difficult to check whether
$\left(RG(\frac{t}{n})E\right)^n$
is uniformly bounded in operator norm with respect to $n$.
Because it is hard to control such iterated applications of
the Gaussian semigroup, one could try to
estimate each factor separately. Then one has to show that
$\left\| RG_0(t)E \right\| \le 1 + ct$ for some $c > 0$, leading
to the upper bound $\e^{ct}$.
The short time diffusion through the boundary, however,
is of order $O(\sqrt{t})$, see~\cite{MPPP07}. This is
why in general only estimates of the kind
$\left\| RG_0(t)E \right\| = 1 + O(\sqrt{t})$ can be obtained.

Almost the same reasoning applies if $\mathrm{C}(\overline{\Omega})$
is replaced by an $L^p$-space, for example by $L^1(\Omega)$.
As $\Omega$ is bounded, uniform convergence already implies
convergence in $L^1(\Omega)$, hence
\[
	T_R(t)u = \lim_{n \to \infty} \bigl( R G_0({\textstyle \frac{t}{n}}) E_\beta \bigr)^n u
	\text{ in } L^1(\Omega)
	\text{ for every } u \in \mathrm{C}(\overline{\Omega})
\]
by what we have already shown. If we want to extend this result
to $u \in L^1(\Omega)$, it suffices to show that the approximating
operators $\bigl( R G_0(\frac{t}{n}) E_\beta \bigr)^n$ remain
bounded in the norm of operators on $L^1(\Omega)$. Here again, there arise
difficulties which are similar to those mentioned in the preceding paragraph
because no non-trivial extension operator from $L^1(\Omega)$ to
$L^1(\mathds{R}^N)$ is contractive.
This shows that for our applications the space $\mathrm{C}(\overline{\Omega})$
has significant advantages.

A related question is whether~\eqref{claim} remains true if the assumption
$\beta \ge 0$ is dropped. We mention that it can be seen that for any $\beta \in L^\infty(\partial \Omega)$
the Laplacian with Robin boundary conditions is the generator of a semigroup on $L^2(\Omega)$
and thus this question makes sense. But if we define $E_\beta$ as in Section~\ref{ExtOp},
we do not even in $\mathrm{C}(\overline{\Omega})$ obtain a contraction if $\beta(z) < 0$
for a point $z \in \partial\Omega$,
causing the same problems again. Moreover, it is clear that
the assumption $\left\| \left( RG_0(t)E_\beta \right)^n \right\|$ which is needed
for Theorem~\ref{chernoff} cannot be fulfilled
since the candidate limit semigroup $T_R(t)$ will not be bounded.
But the latter is merely a problem of rescaling, compare~\cite[Corollary~III.5.3]{Nagel00}.

\bigskip

It should be possible to extend the results to smooth unbounded
open sets $\Omega$ without difficulties because most arguments
are local. However, the other calculations become even more technical.
This is why we have restricted ourselves to bounded domains.

On the other hand, choosing a different (contractive) extension operator
will usually change the situation completely. For example, Theorem~\ref{chernoff}
cannot be applied for the constant extension as in Remark~\ref{Neumann},
reflecting the fact that a worse numerical approximation of the normal derivative
leads to worse convergence behavior. But that extension
operator can be defined even for convex domains without any smoothness assumptions,
which might provide an alternative approximation scheme for less smooth domains.
It is easy to come up with various other extension operators when trying to find
an approximation formula such as~\eqref{claim} for (not necessarily convex)
sets with non-smooth boundary. This is ongoing work and might be the topic
of a future publication.

\bibliographystyle{amsplain}
\bibliography{chernoff_robin}

\end{document}